\documentclass[12pt]{article}

\usepackage{amsxtra,amssymb,amsthm,amsmath,latexsym}
\usepackage{graphicx}
\usepackage{booktabs}
\usepackage{multirow}
\usepackage{afterpage}
\usepackage{setspace}
\usepackage{microtype}
\usepackage{calc,fourier}
\usepackage[sort&compress]{natbib}
\usepackage[T1]{fontenc}
\usepackage[hidelinks]{hyperref}
\usepackage{float}
\usepackage[small]{caption}
\usepackage{tikz}
\usetikzlibrary{positioning}
\allowdisplaybreaks

\numberwithin{equation}{section}

\newcommand{\bee}{\begin{equation*}}
\newcommand{\eee}{\end{equation*}}
\newcommand{\be}{\begin{equation}}
\newcommand{\ee}{\end{equation}}
\newcommand{\ba}{\begin{align}}
\newcommand{\ea}{\end{align}}

\newcommand{\RRR}{\mathbb{R}^3}

\title{A Fast Algorithm for Solving Scalar Wave Scattering Problem by Billions of Particles}
\author{A.G. Ramm\footnote{Mailing address:  Mathematics Department, 138 Cardwell Hall, Manhattan, KS 66506} , N. T. Tran$\dag$ \\
\small Department of Mathematics\\
\small Kansas State University, Manhattan, KS 66506-2602, USA\\
\small \texttt{*ramm@math.ksu.edu} \\
\small \texttt{$\dag$nhantran@math.ksu.edu}
}

\date{}
\begin{document}
\maketitle

\begin{abstract}
Scalar wave scattering by many small particles of arbitrary shapes
with impedance boundary condition is studied. The problem is solved asymptotically and numerically under the assumptions $ a \ll d \ll \lambda$, where $k = 2\pi/\lambda$ is the wave number, $\lambda$ is the wave length, $a$ is the characteristic size of the particles, and $d$ is the smallest distance between neighboring particles. A fast algorithm for solving this wave scattering problem by billions of particles is presented. The algorithm comprises the derivation of the (ORI) linear system and makes use of Conjugate Orthogonal Conjugate Gradient method and Fast Fourier Transform. Numerical solutions of the scalar wave scattering problem with 1, 4, 7, and 10 billions of small impedance particles are achieved for the first time. In these numerical examples, the problem of creating a material with negative refraction coefficient is also described and a recipe for creating materials with a desired refraction coefficient is tested.
\end{abstract}

\noindent\textbf{Key words:} wave scattering; small impedance particles; fast algorithm; billions of particles; negative refraction; meta-materials. \\

\noindent\textbf{MSC:} 35J05; 35J25; 70-08; 68W10; 74Q15.

\section{Introduction} \label{sec1}
Wave scattering is a natural phenomenon that happens in everyday life, for example, light scattering in the atmosphere, light scattering by cosmic dust and by the dust in atmosphere, sound wave scattering by packs of fish in the ocean, etc.
Studying wave scattering is a subject that has attracted much attention from scientists and engineers since it has many practical applications, for example, in  medical image processing,  geophysical prospecting, quantum theory, materials science, etc.
The wave scattering theory gives insights into the structure of the materials, see \cite{R278} and \cite{R470}. Wave scattering by small particles was studied by Lord Rayleigh, see \cite{Ray}, who understood that the main term in the scattered field is the dipole radiation. For particles of an arbitrary shape he did not give formulas for calculating the induced dipole moment with a desired accuracy for bodies of arbitrary shapes. This was done in \cite{R57} and \cite{R476}. Scalar wave scattering by small impedance particles, developed in \cite{R635}, and used in this paper, has practically important physical features:
the field scattered by such particles is $O(a^{2-\kappa})$, as $a\to 0$, which is  much larger than the field in Rayleigh scattering which is $O(a^3)$. Here $a$ is the characteristic size of small particles and it is assumed that the boundary impedance of a particle is $\zeta=ha^{-\kappa}$, where $h$ and $\kappa\in [0,1)$ are constants. The theory of wave scattering by many small impedance particles of an arbitrary shape has been developed in \cite{R632} and \cite{R635}, and is a basis for the computational results in this paper. Our basic physical assumptions are $a<<d<<\lambda$, where $\lambda$ is the wavelength and $d$ is  the minimal distance between neighboring particles. The theory corresponding to the assumptions $a<<\lambda<<d$ is simple and has been used in many cases. It corresponds to the assumption that the effective field in the medium is equal to the incident field. In quantum mechanics it is called the Born approximation, and elsewhere the term weak scattering is used, see \cite{R278} and \cite{R470}.

We do not assume that the particles are distributed in the vertices of a fixed grid with the step size $d$. They can be distributed randomly or not randomly. The small particles can be described by the inequality $ka<<1$, where $k$ is the wave number, $k=2\pi/\lambda$.
In \cite{R635} one can find a detailed presentation of this theory.
In \cite{R632} and \cite{R635} the developed theory has been applied to materials science: it
was proved that by distributing  small particles with  prescribed boundary impedances in a given bounded domain, one can create materials with any desired refraction coefficient, in particular, with negative refraction coefficient, which is of interest for the theory of meta-materials, see \cite{EB2005}.

Earlier numerical results on wave scattering by not more than one million particles, based on the above theory, were reported in \cite{Tran} and \cite{AR2011}. In this paper, for dealing with $10^{10}$ small impedance particles, an essentially novel computational procedure which requires parallel computations at a large scale is developed. The numerical solution of the wave scattering problem with so many small particles, ten billions, is obtained, apparently, for the first time. There are many papers on waves and static fields in the many-body systems. We mention just a few papers \cite{GR1987},\cite{GR1997}, \cite{DardenYorkPedersen}, \cite{PhilWhite}. In these and many other papers in this area the theoretical basis for the computational results is quite different from ours, and there were no computational results on scattering by billions of particles, to our knowledge.  In Section  \ref{sec3} the computational difficulties that we have faced and the methods to overcome these are briefly described.

In  Section \ref{sec2} the theory, on which the computational results are based, is outlined. In  Section \ref{sec3} a fast algorithm for solving wave scattering problem with many small impedance particles of arbitrary shapes is described.
The algorithm is based on 3D convolution, Fast Fourier Transforms (FFT), and Conjugate Orthogonal Conjugate Gradient method (COCG), see \cite{COCG, Clemens} and \cite{CG}. It exploits the structure of the Green's function of the Helmholtz equation in the wave scattering problem and drastically reduces the total number of operations required for solving this problem. The fast computational methods, such as (FFT), have been widely used in various computational
problems, see \cite{HE1988, TDG2004} and \cite{BTK2001}, but the scale of the problem we deal with requires new computational  techniques briefly described in  \ref{sec3}.
Numerical examples are presented in Section \ref{sec4} to illustrate the practical usage of the algorithm. In these numerical examples the algorithm is implemented in parallel and the scalar wave scattering problem is solved with one, four, seven, and ten billions of particles using Gordon super computer at the Extreme Science and Engineering Discovery Environment (XSEDE).

\section{Scalar wave scattering by many small impedance particles} \label{sec2}
Consider a bounded domain $\Omega \subset \RRR$ filled with a material whose refraction coefficient is $n_0(x)$. The assumptions on this coefficient are formulated below \eqref{eq2.5}.
Suppose there are $M$ small particles $D_m$ distributed  in $\Omega$ so that the minimal distance between neighboring particles, $d$, is much greater than the maximal radius of the particles, $a=\frac{1}{2}\max_{1 \leq m \leq M} \text{diam}D_m$, and much less than the wave length, $\lambda$,
$a\ll d \ll \lambda$. Let $D$ be the union of $D_m$, $D := \bigcup_{m=1}^M D_m$, $D\subset \Omega$, and $D':=\RRR \setminus D$ be the exterior domain. Suppose that the boundary impedance of the m\textsuperscript{th} particle is $\zeta_m$, $\zeta_m=\frac{h(x_m)}{a^\kappa}$, where $h(x)$ is a given continuous function in $D$ such that Im $h \leq 0$ in $D$, and $x_m$ is a point inside $D_m$. This point gives the position of the m\textsuperscript{th} particle in $\RRR$. Let $\kappa$ be a given constant, $\kappa \in [0,1)$.
 The scattering problem is formulated as follows:
\begin{align}
    &(\nabla^2+k^2 n_0^2(x))u=0 \quad\text{ in } D', \quad k=\text{const}>0, \quad
     ka \ll 1 \label{eq2.1} \\
    &u_N=\zeta_m u \quad \text{ on }S_m:=\partial D_m, \quad \text{Im }\zeta_m \leq 0,  \quad 1 \leq m \leq M, \label{eq2.2} \\
    &u(x) = u_0(x)+v(x),  \label{eq2.3} \\
    &u_0(x) = e^{ik\alpha \cdot x}, \quad |\alpha|= 1,  \label{eq2.4} \\
    &v_r-ikv=o(1/r), \quad r:=|x| \to \infty. \label{eq2.5}
\end{align}
Here $k$ is the wave number, $k=2\pi/\lambda$, $u_0$ is the incident plane wave, $v$ is the scattered wave, $\alpha$ is the direction of the incident wave, $\vec{N}$ is the outer unit normal to $S_m$, the refraction coefficient  $n_0(x)=1$ in $\Omega':=\RRR \setminus \Omega$. It is assumed that $n_0(x)$ is a Riemann-integrable function and that Im $n_0^2(x) \geq 0$ in $\Omega$. Equation \eqref{eq2.5} is called the radiation condition. It was proved in \cite{R635} that if Im $n_0^2(x) \geq 0$ and Im $h(x) \leq 0$, then the scattering problem  \eqref{eq2.1}-\eqref{eq2.5} has a unique solution and it can be found in the form
\be \label{eq2.6}
    u(x)=u_0(x)+\sum_{m=1}^M \int_{S_m} G(x,t)\sigma_m(t)dt.
\ee
Note that $u$ in \eqref{eq2.6} satisfies \eqref{eq2.1} and \eqref{eq2.5} for any $\sigma_m$. Thus, one just needs to find $\sigma_m$ so that $u$ satisfies \eqref{eq2.2}.

In \eqref{eq2.6}, $G(x,y)$ is the Green's function of the Helmholtz equation \eqref{eq2.1},
$G$ satisfies the equation
\be
	[\nabla^2 + k^2 n_0^2(x)] G = -\delta(x-y) \quad \text{ in } \RRR
\ee
and the radiation condition \eqref{eq2.5}. The functions  $\sigma_m(t)$ are unknown continuous functions.
 These functions are uniquely defined by the boundary condition \eqref{eq2.2}, see \cite{R632}.
If $n_0^2=1$ in $\RRR$, then
\be \label{eq2.6a}
	G(x,y)=\frac{e^{ik|x-y|}}{4\pi|x-y|}.
\ee
The assumption $n_0^2=1$ in $\RRR$ is not a restriction in the problem of creating materials
with a desired refraction coefficient. In the general case, when  $n_0^2$ is a function of $x$,
the Green's function $G$ has to be computed.

From \eqref{eq2.6}, one gets
\be \label{eq2.7}
    u(x)=u_0(x) + \sum_{m=1}^M G(x,x_m)Q_m + \sum_{m=1}^M \int_{S_m} [G(x,t)-G(x,x_m)] \sigma_m(t)dt,
\ee
where
\be \label{eq2.8}
    Q_m:=\int_{S_m} \sigma_m(t)dt.
\ee
It is proved in \cite{R635} that in \eqref{eq2.7}
\be \label{eq2.9}
    |G(x,x_m)Q_m| \gg \left|\int_{S_m} [G(x,t)-G(x,x_m)] \sigma_m(t)dt\right|,
\ee
as $a \to 0$ and $|x-x_m| \geq a$.
Therefore, the solution to the scattering problem can be well approximated by the sum
\be \label{eq2.10}
    u(x) \sim u_0(x)+\sum_{m=1}^M G(x,x_m)Q_m.
\ee
Thus, instead of finding the unknown functions $\sigma_m(t)$ from a system of boundary integral equations, as is usually done when one solves a wave scattering problem,  we just need to find the unknown numbers $Q_m$ to get the accurate approximation of the solution.
This makes it possible to solve problems with so large number of particles that it was not possible to do earlier.

To find the numbers $Q_m$, let us define the effective field $u_e(x)$.
The effective field acting on the j\textsuperscript{th} particle is defined as follows
\be \label{eq2.11}
    u_e(x_j):=u(x_j)-\int_{S_j} G(x_j,t)\sigma_j(t)dt,
\ee
or equivalently
\be \label{eq2.12}
    u_e(x_j)=u_0(x_j)+\sum_{m=1 ,m \neq j}^M \int_{S_m} G(x_j,t)\sigma_m(t)dt,
\ee
where $x_j$ is a point in $D_j$.
The asymptotic formula for $Q_m$ is derived in \cite{R635}:
\be \label{eq2.13}
    Q_m = -c_S a^{2-\kappa} h(x_m)u_e(x_m)[1+o(1)], \quad a \to 0,
\ee
where $c_S>0$ is a constant depending on the shape of the particle,
\be
	|S_m|=c_Sa^2,
\ee
where $|S_m|$ is the surface area of $S_m$. If $S_m$ is a sphere of radius $a$, then $c_S=4\pi$. We assume for simplicity that $c_S$ does not depend on $m$, that is, all the particles are of the same shape.

Let us derive a formula for the effective field. From \eqref{eq2.12}-\eqref{eq2.13} one gets
\be \label{eq2.14}
    u_e(x_j) \simeq u_0(x_j)-c_S\sum_{m=1, m \neq j}^M G(x_j,x_m)h(x_m)u_e(x_m)a^{2-\kappa},
\ee
as $a \to 0$ and $1 \leq j \leq M$.

Denote $u_j:=u_e(x_j),  u_{0j}:=u_0(x_j), G_{jm}:=G(x_j,x_m)$, and $h_m:=h(x_m)$. Then \eqref{eq2.14} can be rewritten as
a linear algebraic system for the unknown numbers $u_m$:
\be \label{eq2.15}
    u_j = u_{0j}-c_S\sum_{m=1, m \neq j}^M G_{jm} h_m a^{2-\kappa} u_m, \quad\text{as } a \to 0, \quad 1 \leq j \leq M.
\ee
In \eqref{eq2.15}, the numbers $u_j$, $1 \le j \le M$, are unknowns. We call \eqref{eq2.15} the  \textit{original linear algebraic system} (ORI). It was proved in \cite{R635} that under the assumptions
\be \label{eq2.16}
    d=O\left(a^{\frac{2-\kappa}{3}}\right), \quad\text{and } M=O\left(\frac{1}{a^{2-\kappa}}\right), \quad\text{for } \kappa \in [0,1),
\ee
the numbers $u_j$, $1 \leq j \leq M,$ can be uniquely found by solving (ORI) for all sufficiently small $a$. If the numbers $u_m$ are known,
then the numbers $Q_m$ can be calculated by formula \eqref{eq2.13} and the approximate solution to the wave scattering problem \eqref{eq2.1}-\eqref{eq2.5} can be computed by \eqref{eq2.10}. {\em This solution is asymptotically exact as $a \to 0$.}

The method for solving many-body wave scattering problem, described above, differs in principle from the Fast Multipole Method (FMM), used in many papers, of which we mention just two: \cite{GR1987} and \cite{GR1997}.
The difference between FMM and our method briefly can be explained as follows: the theoretical basis is different, our method is developed
for scattering by small impedance particles of arbitrary shapes and is based on the asymptotically exact formula for the field, scattered by one small particle, and on the assumption $d \gg a$; and we derive an integral equation for the limiting field in the medium consisting of many small particles as $a\to 0$. We do not use multipole expansions. One of the advantages of our method is in the asymptotic exactness of the method as $a \to 0$.

Next, let us derive the reduced order linear system for solving the wave scattering problem. Let $\Delta$ be an arbitrary subdomain of $\Omega$. Assume that the distribution of particles in $\Delta$ satisfies this law
\be \label{eq2.17}
    \mathcal{N}(\Delta)=\frac{1}{a^{2-\kappa}} \int_{\Delta} N(x)dx[1+o(1)], \quad\text{as } a \to 0.
\ee
Here $N(x) \ge 0$ is a given continuous function in $\Omega$. The function  $N(x)$ and the number $\kappa\in [0,1)$ can be chosen
by the experimenter as he(she) desired. The number $\mathcal{N}(\Delta)$ is the total number of the embedded particles in $\Delta$.

Let $\Omega$ be partitioned into $P$ non-intersecting sub-cubes $\Delta_p$ of side $b$ such that $b \gg d \gg a$, where $b=b(a)$, $d=d(a)$, and $\lim_{a\to 0}\frac{d(a)}{b(a)}=0$. Here $P \ll M$,  and each sub-cube contains many particles.  If the function $N(x)$  in \eqref{eq2.17} is continuous and $b \ll 1$, then
\be \label{eq2.18}
	\mathcal{N}(\Delta_p) a^{2-\kappa}=N(x_p)|\Delta_p|[1+o(1)]=a^{2-\kappa}\sum_{x_m \in \Delta_p}1, \quad\text{as } a \to 0,
\ee
where $|\Delta_p|$ is the volume of $\Delta_p$ and $x_p\in \Delta_p$ is an arbitrary point, for example, the center of $\Delta_p$.
Thus, \eqref{eq2.15} can be rewritten as
\be \label{eq2.19}
    u_q = u_{0q}-c_S\sum_{p=1, p \neq q}^P G_{qp} h_p N_p u_p |\Delta_p|, \quad\text{for } 1 \leq q \leq P,
\ee
where $N_p:=N(x_p)$ and $x_p$ is a point in $\Delta_p$, for example, the center of $\Delta_p$. We call \eqref{eq2.19} the \textit{reduced linear algebraic system} (RED). This system is much easier to solve since $P \ll M$.

Let $|\Delta_p| \to 0$. Then it follows from \eqref{eq2.19} that the limiting integral equation for $u=u(x)$ holds
\be \label{eq2.20}
    u(x)=u_0(x)-c_S\int_\Omega G(x,y)h(y)N(y)u(y)dy, \quad\text{for } x \in \mathbb{R}^3,
\ee
 if the assumption \eqref{eq2.17} is satisfied. The sum in \eqref{eq2.19} is the Riemannian sum for the integral in \eqref{eq2.20} which converges to this integral when $\max_{p}|\Delta_p|\to 0$ (see \cite{R635} for the proof of convergence).

Let
\be \label{eq2.20a}
	p(x):=c_SN(x)h(x).
\ee
Then \eqref{eq2.20} can be written as
\be \label{eq2.21}
    u(x)=u_0(x)-\int_\Omega G(x,y)p(y)u(y)dy, \quad\text{for } x \in \mathbb{R}^3.
\ee
 Here $u=u(x)$ is the limiting field in the medium created by embedding many small impedance particles distributed according to
  equation \eqref{eq2.17}. We call \eqref{eq2.21} the \textit{limiting integral equation} (IE).

Now, applying the operator $(\nabla^2+k^2n_0^2)$ to \eqref{eq2.21} and using the equation $(\nabla^2+k^2n_0^2)G(x,y)=-\delta(x-y)$, one gets
\be  \label{eq2.22}
    (\nabla^2+k^2n_0^2)u(x)=p(x)u(x).
\ee
This implies
\be \label{eq2.23}
	(\nabla^2+k^2n^2)u=0,
\ee
where
\be \label{eq2.24}
	n^2(x):=n_0^2(x)-k^{-2}p(x),
\ee
and {\em  $n(x)$ is the new refraction coefficient of the limiting medium.}
 Since Im$h(x) \le 0$ and Im$n_0^2(x) \ge 0$, one concludes that Im$n^2(x) \ge 0$. From \eqref{eq2.24}, one gets
\be \label{eq2.25}
	p(x)=k^2[n_0^2(x)-n^2(x)].
\ee
By equation \eqref{eq2.20a}, $h(x)$ can be computed as
\be \label{eq2.26}
	h(x)=\frac{p(x)}{c_SN(x)}.
\ee
{\em This gives  a method for creating new materials with a desired refraction coefficient $n(x)$ by embedding many small impedance particles into a given material with the original refraction coefficient $n_0$ using the distribution law \eqref{eq2.17}.}

\section{A fast algorithm for solving wave scattering problem by billions of particles} \label{sec3}
For solving the (RED) linear system, one can use any iterative method, namely GMRES, see \cite{GMRES}. Since the order of (RED) can be made much smaller than that of (ORI), the computation is very fast. Therefore, our remaining goal is to develop a fast algorithm for solving (ORI) in order to get the solution of the scattering problem \eqref{eq2.1}-\eqref{eq2.5}. Our algorithm is a combination of the Conjugate Orthogonal Conjugate Gradient (COCG) method, 3D convolution, and FFT. When one solves a linear algebraic system using iterations, matrix-vector multiplications are carried out in the iterative process. These multiplications take most of the computation time. If the linear system is very large, it takes a huge amount of time to finish only one matrix-vector multiplication in a standard way, since this multiplication is of the order $O(n^2)$. In some cases it is practically impossible to perform such computations, for example, when the system is dense and has more than one billion equations and unknowns. In this section we present an algorithm that  greatly reduces the total number of operations (from $O(n^2)$ to $O(n\log n)$) and decreases the overall computation time of the iterative process by handling the matrix-vector multiplication by using 3D convolution and FFT.

There are numerous methods which also employ FFT to solve different problems, for example, Precorrected-FFT method for electrostatic analysis of complicated 3D structures \cite{PhilWhite}, or Particle mesh Ewald method for Ewald sums in large systems \cite{DardenYorkPedersen}, etc. Nevertheless, none of these papers deals with the scale that we face solving the scalar wave scattering problem with ten billion particles, i.e., solving a $10^{10} \times 10^{10}$ linear system. This is done for the first time in our work.  We have to develop a new algorithm that can solve two major problems in our computing: memory and time. First, it is impossible to store a $10^{10} \times 10^{10}$ dense matrix in any currently available super computer. Suppose we use only single precision. Then it would take 800 million terabytes of memory to store only one matrix to do the computation, since each complex number is 8 bytes. Furthermore, we will suffer network latency and traffic jams which cause a halt in our computation if we use such an amount of memory in any parallel cluster. Second, it is impossible to do the computation at order $O(n^2)$ for a $10^{10} \times 10^{10}$ linear system in a reasonable and permitted time. We deal with the first computing problem, memory, by finding a way to store the $n \times n$ matrix in a 3D cube which is equivalent to only one $n \times 1$ vector in size and avoid moving terabytes of data around. The second computing problem, time, is resolved by reducing the number of operations from $O(n^2)$ to $O(n\log n)$. The details on how to deal with these two major computing problems at our scale and how to set up the wave scattering problem in order to solve it in parallel clusters are described in this section.

Consider the following summation in the original linear system \eqref{eq2.15} of the wave scattering problem
\be \label{eq3.3.1}
	\sum_{m=1, m \ne j}^M G(\bold{x_j},\bold{x_m})u(\bold{x_m}),
\ee
where by bold letters vectors are denoted. $G$ in equation \eqref{eq3.3.1} is the Green's function of the form
\be \label{eq3.3.2}
	G(\bold{x_j},\bold{x_m})=\frac{e^{ik|\bold{x_j}-\bold{x_m}|}}{4\pi|\bold{x_j}-\bold{x_m}|},
\ee
where $\bold{x_j}, \bold{x_m}$ are the positions of the j\textsuperscript{th} and m\textsuperscript{th} particles in $\RRR$, respectively. If we write $G(\bold{x}-\bold{y}):=G(\bold{x},\bold{y})$, the summation in \eqref{eq3.3.1} will be
\be \label{eq3.3.4}
	\sum_{m=1, m \ne j}^M G(\bold{x_j}-\bold{x_m})u(\bold{x_m}),
\ee
which is a discrete convolution of $G$ and $u$, $G*u$, if $m \ne j$ is dropped.

In the linear system \eqref{eq2.15}, $G$ is an $M\times M$ matrix, where $M$ is the total number of particles, and $u$ is an $M \times 1$ vector. When solving the linear system \eqref{eq2.15} using COCG iterative algorithm, the matrix-vector multiplication in \eqref{eq3.3.4} needs to be executed. If we do this matrix-vector multiplication in the standard way, it would take $O(M^2)$ operations. This is very expensive in terms of computation time if $M$ is very large, for example,  $M \ge 10^{6}$. Therefore, we have to find a new way to do the matrix-vector multiplication.

The convolution in \eqref{eq3.3.4} can be carried out by using Convolution theorem as follows:
\be \label{eq3.3.8}
	G*u=\mathcal{F}^{-1}(\mathcal{F}(G*u))=\mathcal{F}^{-1}(\mathcal{F}(G)\cdot \mathcal{F}(u)),
\ee
where the $\cdot$ stands for the component-wise multiplication of two vectors and its result is a vector.

If particles are distributed uniformly, one can just use FFT to quickly compute this convolution. Otherwise, one can use Nonequispaced Fast Fourier Transform (NFFT), see \cite{NFFT}, \cite{KKP2009}, and \cite{PP2013}. Our method is valid for variable $N(x)$, see formulas \eqref{eq2.19} and \eqref{eq2.20}. Alternatively, one can just use the (RED) linear system or (IE) with much lower order to solve the wave scattering problem without using FFT nor NFFT. The theory in section \ref{sec2} shows that the solution to (RED) or (IE) yields a solution to (ORI) with high accuracy, the error tends to zero as $a \to 0$.

To illustrate the idea, let us assume for simplicity that particles are distributed uniformly, that is, $N(x)=$const in \eqref{eq2.17}.
Let $\bold{m}=(x_m,y_m,z_m)$ be the position of the m\textsuperscript{th} particle in $\RRR$, where $x_m$, $y_m$, and $z_m$ are real numbers. We will assume our domain is a unit cube (different domains can be treated similarly),  this cube is placed in the first octant and the origin is one of its vertices, then $\bold{m}$ can be rewritten as a product of the scalar factor $d>0$ and a vector $(m_1, m_2, m_3)$:
\be \label{eq3.3.5}
	\bold{m}=d(m_1,m_2,m_3),
\ee
where $d$ is the distance between neighboring particles, a scalar,  $(m_1,m_2,m_3)$ is a vector whose components $m_1, m_2$, and $m_3$ are integers in $[0,b)$, and $b=M^{1/3}$ is the number of particles on a side of the cube.

In the convolution \eqref{eq3.3.4} suppose that $\bold{x_j}=d(j_1,j_2,j_3)$ and $\bold{x_m}=d(m_1,m_2,m_3)$, one can write \eqref{eq3.3.4} as
\begin{align}
	G*u &= \sum_{m=1, m \ne j}^M  G(\bold{x_j}-\bold{x_m})u(\bold{x_m})  \label{eq3.3.6} \\
		&= \sum_{{\footnotesize \begin{array}{c}
			m_1,m_2,m_3=0 \\
			(m_1,m_2,m_3) \ne (j_1,j_2,j_3)
		\end{array}}}^{b-1}
		G(j_1-m_1,j_2-m_2, j_3-m_3)u(m_1,m_2,m_3).  \label{eq3.3.7}
\end{align}
This is a 3D convolution of $G$ and $u$.

\begin{figure}[htbp]
\centering
\includegraphics[width=0.3\linewidth]{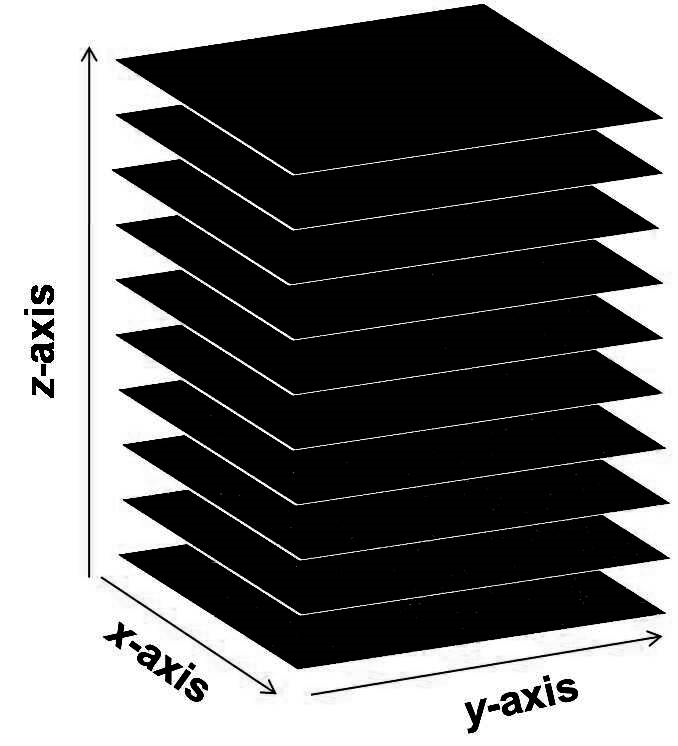}
\caption{Dividing a cube into a stack of planes for storing across machines.}
\label{fig0}
\end{figure}
In order to do this convolution, we need to store matrix $G$ as a vector. This reduces drastically the amount of memory for storing the original $M\times M$ matrix to a much smaller amount for storing an $M \times 1$ vector, which is also denoted by $G$. Since this vector depends on the three components, $G=G(j_1-m_1,j_2-m_2, j_3-m_3)$ where $j_1, m_1, j_2, m_2, j_3, m_3$ are integers in $[0,b)$ and $(j_1,j_2,j_3) \ne (m_1,m_2,m_3)$, we can alternatively store it as a cube of size $b\times b\times b$. Similarly, vector $u(m_1,m_2,m_3)$ is also stored as a cube of size $b\times b\times b$. Each cube is a stack of planes which will be distributed across all machines in a cluster for parallelizing the computations and reducing moving data around to prevent traffic jam and network latency, see Figure \ref{fig0}. Each machine will work on its local data. Information can be shared among machines but as minimal as possible. Matrix-vector multiplication is done via a function handle without storing any matrix.

When implementing the 3D convolution we need to pad the cubes $G$ and $u$ as follows:
\begin{itemize}
	\item Pad in x-direction
	\item Pad in y-direction
	\item Pad in z-direction
\end{itemize}
\begin{center}
\begin{tikzpicture}[on grid,scale=0.3]
  \shade[yslant=-0.5,right color=gray!10, left color=black!50]
    (0,0) rectangle +(3,3);
  \draw[yslant=-0.5] (0,0) grid (3,3);
  \shade[yslant=0.5,right color=gray!70,left color=gray!10]
    (3,-3) rectangle +(3,3);
  \draw[yslant=0.5] (3,-3) grid (6,0);
  \shade[yslant=0.5,xslant=-1,bottom color=gray!10,
    top color=black!80] (6,3) rectangle +(-3,-3);
  \draw[yslant=0.5,xslant=-1] (3,0) grid (6,3);
\end{tikzpicture}
$\xrightarrow{Pad \; x}$
\begin{tikzpicture}[on grid,scale=0.3]
  \shade[yslant=-0.5,right color=gray!10, left color=black!50]
    (1,1) rectangle +(3,3);
  \draw[yslant=-0.5] (0,0) grid (3,3);
  \shade[yslant=0.5,right color=gray!70,left color=gray!10]
    (4,-3) rectangle +(3,3);
  \draw[yslant=0.5] (3,-3) grid (7,0);
  \shade[yslant=0.5,xslant=-1,bottom color=gray!10,
    top color=black!80] (7,3) rectangle +(-3,-3);
  \draw[yslant=0.5,xslant=-1] (3,0) grid (7,3);
\end{tikzpicture}
$\xrightarrow{Pad \; y}$
\begin{tikzpicture}[on grid,scale=0.3]
  \shade[yslant=-0.5,right color=gray!10, left color=black!50]
    (1,1) rectangle +(3,3);
  \draw[yslant=-0.5] (0,0) grid (4,3);
  \shade[yslant=0.5,right color=gray!70,left color=gray!10]
    (4,-3) rectangle +(3,3);
  \draw[yslant=0.5] (3,-3) grid (7,0);
  \shade[yslant=0.5,xslant=-1,bottom color=gray!10,
    top color=black!80] (7,3) rectangle +(-3,-3);
  \draw[yslant=0.5,xslant=-1] (3,-1) grid (7,3);
  \draw[yslant=0.5] (4,-4) grid (8,-1);
\end{tikzpicture}
$\xrightarrow{Pad \; z}$
\begin{tikzpicture}[on grid,scale=0.3]
  \shade[yslant=-0.5,right color=gray!10, left color=black!50]
    (1,1) rectangle +(3,3);
  \draw[yslant=-0.5] (0,0) grid (4,4);
  \shade[yslant=0.5,right color=gray!70,left color=gray!10]
    (4,-3) rectangle +(3,3);
  \draw[yslant=0.5] (3,-3) grid (7,0);
  \shade[yslant=0.5,xslant=-1,bottom color=gray!10,
    top color=black!80] (7,3) rectangle +(-3,-3);
  \draw[yslant=0.5,xslant=-1] (3,-1) grid (7,3);
  \draw[yslant=0.5] (4,-4) grid (8,0);
  \draw[yslant=0.5,xslant=-1] (4,0) grid (8,4);
\end{tikzpicture}
\end{center}
For each direction, the padding is illustrated by this example \\
For $G$:
\begin{center}
\boxed{
	\begin{array}{cccc}
		1 & 2 & 3 & 4 \\
		5 & 6 & 7 & 8 \\
		9 & 10 & 11 & 12 \\
		13 & 14 & 15 & 16
	\end{array}
}
\quad$\xrightarrow{Pad}$\quad
\boxed{
	\begin{array}{cccccc}
		1 & 2 & 3 & 4 & 3 & 2 \\
		5 & 6 & 7 & 8 & 7 & 6 \\
		9 & 10 & 11 & 12 & 11 & 10 \\
		13 & 14 & 15 & 16 & 15 & 14
	\end{array}
}
\end{center}
For $u$:
\begin{center}
\boxed{
	\begin{array}{cccc}
		1 & 2 & 3 & 4 \\
		5 & 6 & 7 & 8 \\
		9 & 10 & 11 & 12 \\
		13 & 14 & 15 & 16
	\end{array}
}
\quad$\xrightarrow{Pad}$\quad
\boxed{
	\begin{array}{cccccc}
		1 & 2 & 3 & 4 & 0 & 0 \\
		5 & 6 & 7 & 8 & 0 & 0 \\
		9 & 10 & 11 & 12 & 0 & 0 \\
		13 & 14 & 15 & 16 & 0 & 0
	\end{array}
}
\end{center}
This means that we pad $G$ using its entries and pad $u$ with zeros. As described in the example above, for padding $G$ we copy all columns except the first and the last ones and put them symmetrically through the last column. This will create a periodic signal $G$. In fact, if one places padded $G$ continuously, one can see a periodic signal. Since we only need to perform linear convolution on $M$-length vectors, the result we need is an $M$-length vector. All the entries after the $M$-th entry in the convolution will be discarded. So, we pad $u$ with zeros just to have the same length with the padded $G$ to do the cyclic convolution in computer. Cyclic convolutions allow us to compute linear convolutions by means of Discrete Fourier Transforms (DFT). After padding $G$ and $u$ will have size $(2b-2)^3$.

The Fourier transform and inverse Fourier transform are of order $O(n\log n)$, and vector pointwise multiplication is of order $O(n)$ if the vectors are $n \times 1$. In our case the total number of operations for computing $G*u=\mathcal{F}^{-1}(\mathcal{F}(G)\cdot \mathcal{F}(u))$ is
\be \label{eq3.3.9}
	n\log n+n\log n+n+n\log n=O(n\log n), \quad n=(2b-2)^3,
\ee
since the Fourier transforms $\mathcal{F}(G)$ and $\mathcal{F}(u)$ are of order $O(n\log n)$, the vector point-wise multiplication  $\mathcal{F}(G)\cdot \mathcal{F}(u)$ is of order $O(n)$, and the inverse Fourier transform $\mathcal{F}^{-1}(\mathcal{F}(G)\cdot \mathcal{F}(u))$ is of order $O(n\log n)$.
If we compare this with the standard matrix-vector multiplication which takes $M^2$ operations $(M=b^3)$, this is a huge reduction of the number of operations and computation time, when $M$ is very large, say $M \ge 10^{9}$.

This algorithm is applicable not only to solving scalar wave scattering problems but also to other PDE problems, for example, aeroacoustics, signal processing, propagator in quantum mechanics and quantum field theory, etc.

\section{Numerical examples} \label{sec4}
The algorithm described in Section \ref{sec3} is implemented in parallel using the Portable, Extensible Toolkit for Scientific Computation (PETSc) library developed at Argonne National Laboratory (ANL), see \cite{Petsc}. For implementing FFT, Fastest Fourier Transform in the West (FFTW) library is used, see \cite{FFTW}. The wave scattering problem is solved using Gordon super computer at XSEDE. "Gordon is a dedicated XSEDE cluster designed by Appro and SDSC consisting of 1024 compute nodes and 64 I/O nodes. Each compute node contains two 8-core 2.6 GHz Intel EM64T Xeon E5 (Sandy Bridge) processors and 64 GB of DDR3-1333 memory", see \cite{Gordon}. Table \ref{tab1} shows the technical information of one compute node in Gordon.
\begin{table}[htb]
  \centering
  \caption{Compute node Intel EM64T Xeon E5.}
    \begin{tabular}{rr}
    \toprule
    \multicolumn{1}{c}{System Component} & \multicolumn{1}{c}{Configuration} \\
    \midrule
    \multicolumn{1}{l}{Sockets} & \multicolumn{1}{l}{2} \\
    \multicolumn{1}{l}{Cores} & \multicolumn{1}{l}{16} \\
    \multicolumn{1}{l}{Clock speed} & \multicolumn{1}{l}{2.6 GHz} \\
    \multicolumn{1}{l}{Flop speed} & \multicolumn{1}{l}{333 Gflop/s} \\
    \multicolumn{1}{l}{Memory capacity} & \multicolumn{1}{l}{64 GB DDR3-1333} \\
    \multicolumn{1}{l}{Memory bandwidth} & \multicolumn{1}{l}{85 GB/s} \\
    \multicolumn{1}{l}{STREAM Triad bandwidth} & \multicolumn{1}{l}{60 GB/s} \\
    \bottomrule
    \end{tabular}%
  \label{tab1}%
\end{table}%
"The network topology of Gordon is a 4x4x4 3D torus with adjacent switches connected by three 4x QDR InfiniBand links (120 Gbit/s). Compute nodes (16 per switch) and I/O nodes (1 per switch) are connected to the switches by 4x QDR (40 Gbit/s). The theoretical peak performance of Gordon is 341 TFlop/s", see \cite{Gordon}. Table \ref{tab2} shows information about the network of Gordon.
\begin{table}[htb]
  \centering
  \caption{Network summary.}
    \begin{tabular}{rr}
    \toprule
    \multicolumn{2}{c}{QDR InfiniBand Interconnect} \\
    \midrule
    \multicolumn{1}{l}{Topology} & \multicolumn{1}{l}{3D Torus} \\
    \multicolumn{1}{l}{Link bandwidth} & \multicolumn{1}{l}{8 GB/s (bidirectional)} \\
    \multicolumn{1}{l}{MPI latency} & \multicolumn{1}{l}{1.3 $\mu$s} \\
    \bottomrule
    \end{tabular}%
  \label{tab2}%
\end{table}%

The program code is written in C \& C++, compiled with Intel compiler, and linked with MPI library MVAPICH2. The code uses 64-bit integers and single precision. The relative error tolerance used for the convergence of COCG iterations is $2\times 10^{-5}$.

We assume that the domain $\Omega$, which contains all the particles, is a unit cube, placed in the first octant such that the origin is one of its vertices, and particles are distributed uniformly in $\Omega$. Suppose we want to create a new meta-material with the refraction coefficient $n(x)=-1$ in $\Omega$  given a material with the refraction coefficient $n_0(x)=1$ by embedding many small particles into the given material. We assume the particles are spheres, so $c_S=4\pi$. The new refraction coefficient is computed by the following formula
\be \label{eq4.1}
	n(x)=[n_0^2(x)-k^{-2}c_Sh(x)N(x)]^{1/2},
\ee
where $N(x)$ and $h(x)$ are at our choices. For simplicity we choose $N(x)=1$. The choice of $h(x)$ is subject to the physical condition
Im$h\le 0$. If Im$h(x) \le 0$ and Im$n_0^2\ge 0$, then Im$n^2(x) \ge 0$.
The square root in formula \eqref{eq4.1} is of the form
\be \label{eq4.2}
	z^{1/2}=|z|^{1/2}e^{i\frac{\phi}{2}}, \quad \phi :=\arg z, \quad \phi \in [0,2\pi].
\ee
Formula \eqref{eq4.2} defines a one-valued branch of analytic function $z^{1/2}$ in the complex plane  with the cut $[0, +\infty)$.
If one wants to get $n=Be^{i(\pi -\epsilon)}$, where $B>0$ and $\epsilon>0$, then $n^2=B^2e^{i(2\pi-2\epsilon)}$.
When $\epsilon >0$ is very small, one gets practically negative refraction coefficient $n$.
In this experiment, we choose Im$n=0.001$. This violates the assumption  Im$h(x) \le 0$. To justify this violation for very small
values of  Im$h(x)$ we argue as follows.
The integral equation \eqref{eq2.20} is an equation with compact integral operator $T$
\be
	Tu:=c_S\int_D G(x,y)h(y)N(y)u(y)dy.
\ee
It is of Fredholm type with index zero. It is proved in \cite{R635} that equation \eqref{eq2.20}
has at most one solution for Im$h\le 0$. Therefore, the inverse operator $(I+T)^{-1}$ is bounded for  Im$h\le 0$. The set of boundedly invertible operators is open. Therefore the inverse operator $(I+T)^{-1}$
exists and is bounded also for sufficiently small Im$h\ge 0$.

The radius $a$ of the particles and the distance $d$ between neighboring particles are chosen so that
\be \label{eq4.3}
    d=\frac{1}{M^{1/3}}=a^{\frac{2-\kappa}{3}}, \quad\text{and}\quad M=\frac{1}{a^{2-\kappa}},
\ee
where $M$ is the total number of particles embedded in the domain $\Omega$. To solve (IE), we use  a collocation method, dividing the domain into many sub-cubes, taking the collocation points as the centers of these cubes, and then approximating the integral equation by the corresponding Riemannian sum.

The following physical parameters are used to conduct the experiment:
\begin{itemize}
	\item Speed of wave, $v=$ 34400 cm/sec,
	\item Frequency, $f= 1000$ Hz,
	\item Wave number, $k= 0.182651$ cm$^{-1}$,
	\item Direction of plane wave, $\alpha= (1, 0, 0)$,
	\item The constant $\kappa= 0.5$,
	\item Volume of the domain that contains all particles, $|\Omega|= 1$ cm$^3$,
	\item Distribution of particles, $N=Ma^{2-\kappa}/|\Omega|= 1$, i.e. particles are distributed uniformly in the unit cube,
	\item Function $h(x) =$ 2.65481E-09 + i5.30961E-06,
	\item Original refraction coefficient, $n_0=$ 1+i0,
	\item Desired refraction coefficient, $n =$ -1+i0.001\\
	($n=[n_0^2(x)-k^{-2}c_Sh(x)N(x)]^{1/2} = [1^2-0.182651^{-2}4\pi$(2.65481E-09 + i5.30961E-06)1]$^{1/2} =$ -1+i0.001),
	\item Number of small subcubes after partitioning the domain $\Omega$ for solving (RED), $P=8000$.
	\item Number of collocation points for solving (IE), $C= 64000$.
\end{itemize}

\begin{figure}[htb]
\centering
\includegraphics[width=0.7\linewidth]{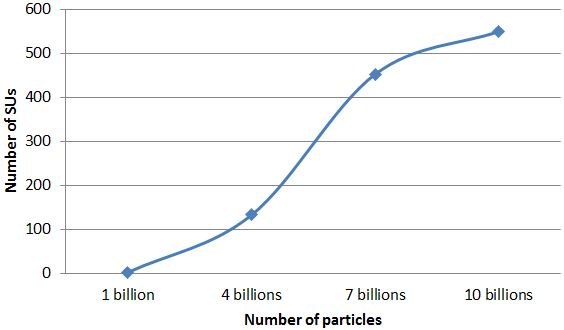}
\caption{Time usage for solving the wave scattering problem.}
\label{fig1}
\end{figure}

\begin{figure}[htbp]
\centering
\includegraphics[width=0.7\linewidth]{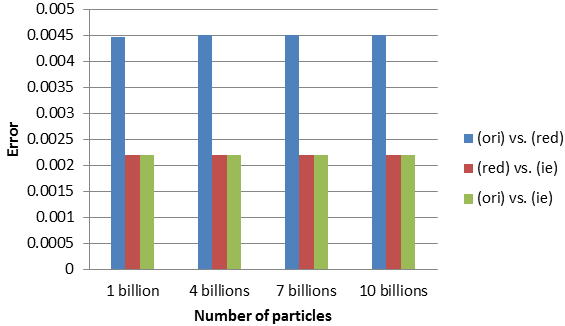}
\caption{Solution differences (errors).}
\label{fig2}
\end{figure}

Table \ref{tab3} and Figure \ref{fig1} show the time usage in Gordon for solving the wave scattering problem with 1, 4, 7, and 10 billion particles using the algorithm described in Section \ref{sec3}. The computation time is measured by Service Unit (SU), 1 SU corresponds to 1 hour/core.
\begin{table}[htb]
  \centering
  \caption{Time usage for solving the wave scattering problem.}
    \begin{tabular}{lcccc}
    \toprule
    Number of particles & 1 billion & 4 billions & 7 billions & 10 billions \\
    \midrule
    Time usage (second) & 103   & 1076  & 2082  & 1674 \\
    Node usage & 8     & 28    & 49    & 74 \\
    Number of SUs & 3.66 & 133.90 & 453.41 & 550.56 \\
    \bottomrule
    \end{tabular}%
  \label{tab3}%
\end{table}%
Table \ref{tab4} and Figure \ref{fig2} show the differences (errors) between the solutions of (ORI) vs. (RED), (RED) vs. (IE), and (ORI) vs. (IE). Since the numbers of unknowns in (ORI), (RED), and (IE) are different, $M \gg P$ and $P < C$, we use interpolation procedure to compare their solutions. For example, let $x$ and $y$ be the solutions of (ORI) and (RED), respectively. We find all the particles $x_i$ that lie in the subcube $\Delta_q$ corresponding to $y_q$ and then find the difference $|x_i-y_q|$. Then, we compute
\be
	\sup_{y_q} \frac{1}{\mathcal{N}(\Delta_q)}\sum_{x_i \in \Delta_q}|x_i-y_q|,
\ee
where $\mathcal{N}(\Delta_q)$ is the number of particles in the subcube $\Delta_q$. This gives the difference between the solutions of (ORI) and (RED). The solution differences between (RED) vs. (IE) and (ORI) vs. (IE) are computed similarly. The numbers in table \ref{tab4} are rounded to the nearest ten-thousandths.
\begin{table}[htb]
  \centering
  \caption{Solution differences (errors).}
    \begin{tabular}{lcccc}
    \toprule
    Number of particles & 1 billion & 4 billions & 7 billions & 10 billions \\
    \midrule
    (ORI) vs. (RED) & 0.0045 & 0.0045 & 0.0045 & 0.0045 \\
    (RED) vs. (IE) & 0.0022 & 0.0022 & 0.0022 & 0.0022 \\
    (ORI) vs. (IE) & 0.0022 & 0.0022 & 0.0022 & 0.0022 \\
    \bottomrule
    \end{tabular}%
  \label{tab4}%
\end{table}%

\begin{figure}[H]
	\centering
	\includegraphics[width=1.05\linewidth,height=6.8cm]{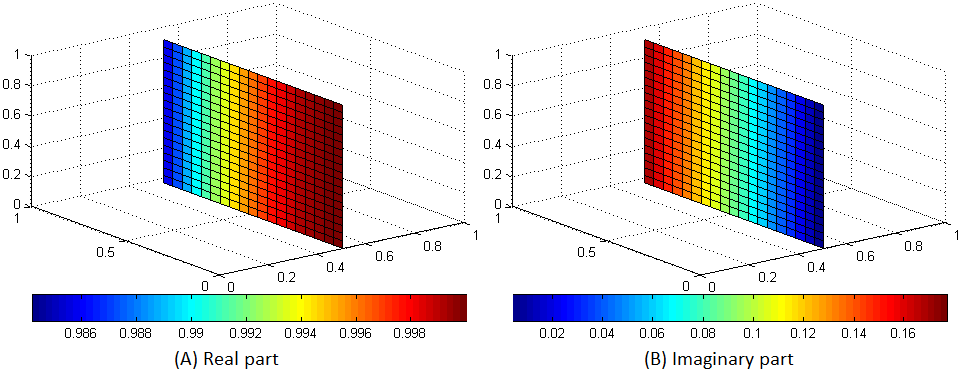}
	\caption{Solution of (ORI) when $M=10^9$.}
	\label{fig3}
\end{figure}

\begin{figure}[H]
	\centering
	\includegraphics[width=1.05\linewidth,height=6.8cm]{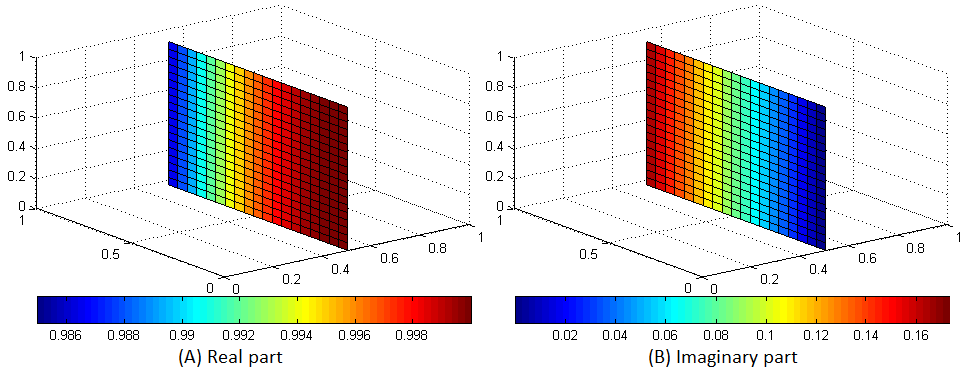}
	\caption{Solution of (RED) when $M=10^9$ and $P=8000$.}
	\label{fig4}
\end{figure}

\begin{figure}[H]
	\centering
	\includegraphics[width=1.05\linewidth,height=6.8cm]{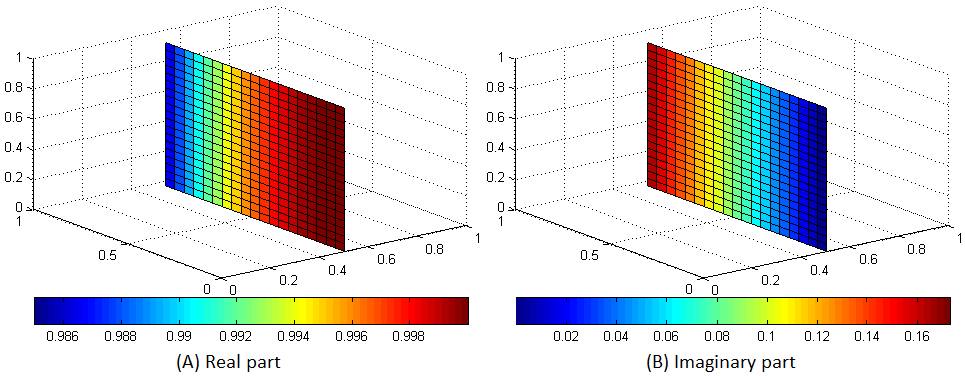}
	\caption{Solution of (IE) when $M=10^9$, $P=8000$, and $C=64000$.}
	\label{fig5}
\end{figure}

For example, Figures \ref{fig3}, \ref{fig4} and \ref{fig5} display vertical slice planes of the solutions, scattering fields, of (ORI), (RED), and (IE), respectively, at the center of the domain $\Omega$, when $M=10^9$ particles, $P=8000$ subcubes, and $C=64000$ collocation points. The relative errors of the convergence of the solutions to (ORI), (RED), and (IE) are 1.72448E-05, 1.94613E-05, and 1.93914E-05, respectively. The relative errors are residual-based and computed using this ratio $\frac{||r||}{||RHS||}$, where $r$ is the residual vector and $RHS$ is the right-hand-side vector. The solution differences between (ORI) vs. (RED), (RED) vs. (IE), and (ORI) vs. (IE) are 0.0045, 0.0022, and 0.0022, respectively. The color bars indicate the  values of the corresponding colors. The values used here are the real part and imaginary part of the scattering fields at the grid points on the slices. 

For reference, tables \ref{tab5}, \ref{tab6}, and \ref{tab7} show the solutions of (ORI), (RED), and (IE), respectively, at the grid points $5\times 5\times 5$ in the unit cube $\Omega$.

\section{Conclusions}
The numerical results in this paper allow one to solve (ORI) for $1\le M \le 10^{10}$. These results show that the solution by (RED) for $M=10^{10}$ agrees with the solution by (ORI) with high accuracy (99.55 \%), and agrees with the solution of (IE) also with high accuracy (99.78\%). Therefore, practically for solving problems with $M \ge 10^6$ one may use (RED) or (IE). For solving the scattering problem for $M<10^6$ numerically one can use (ORI). The accuracy of our numerical method is high if the quantity $ka+ad^{-1}$ is small. Furthermore, it is important to note that this method is of the same order as that of Fast Multipole Method (FMM), $O(n\log n)$. However, we do not use multipole expansions and our method is relatively easy to implement compared to FMM and it can give an asymptotically exact solution to the wave scattering problem.

\pagebreak

\begin{table}[H]
  \centering
  \caption{Solution of (ORI) at the grid points 5$\times$5$\times$5 in the cube.}
  {\fontsize{8.8}{9.6}\selectfont
    \begin{tabular}{rrrrr}
    \toprule
    0.999990+0.004392i  &  0.999990+0.004392i  &  0.999990+0.004392i  &  0.999990+0.004392i  &  0.999990+0.004392i  \\
    0.999990+0.004392i  &  0.999990+0.004392i  &  0.999990+0.004392i  &  0.999990+0.004392i  &  0.999990+0.004392i  \\
    0.999990+0.004392i  &  0.999990+0.004392i  &  0.999990+0.004392i  &  0.999990+0.004392i  &  0.999990+0.004392i  \\
    0.999990+0.004392i  &  0.999990+0.004392i  &  0.999990+0.004392i  &  0.999990+0.004392i  &  0.999990+0.004392i  \\
    0.999990+0.004392i  &  0.999990+0.004392i  &  0.999990+0.004392i  &  0.999990+0.004392i  &  0.999990+0.004392i  \\
          &       &       &       &  \\
    0.999163+0.040911i  &  0.999163+0.040911i  &  0.999163+0.040911i  &  0.999163+0.040911i  &  0.999163+0.040911i  \\
    0.999163+0.040911i  &  0.999163+0.040911i  &  0.999163+0.040911i  &  0.999163+0.040911i  &  0.999163+0.040911i  \\
    0.999163+0.040911i  &  0.999163+0.040911i  &  0.999163+0.040911i  &  0.999163+0.040911i  &  0.999163+0.040911i  \\
    0.999163+0.040911i  &  0.999163+0.040911i  &  0.999163+0.040911i  &  0.999163+0.040911i  &  0.999163+0.040911i  \\
    0.999163+0.040911i  &  0.999163+0.040911i  &  0.999163+0.040911i  &  0.999163+0.040911i  &  0.999163+0.040911i  \\
          &       &       &       &  \\
    0.997002+0.077375i  &  0.997002+0.077375i  &  0.997002+0.077375i  &  0.997002+0.077375i  &  0.997002+0.077375i  \\
    0.997002+0.077375i  &  0.997002+0.077375i  &  0.997002+0.077375i  &  0.997002+0.077375i  &  0.997002+0.077375i  \\
    0.997002+0.077375i  &  0.997002+0.077375i  &  0.997002+0.077375i  &  0.997002+0.077375i  &  0.997002+0.077375i  \\
    0.997002+0.077375i  &  0.997002+0.077375i  &  0.997002+0.077375i  &  0.997002+0.077375i  &  0.997002+0.077375i  \\
    0.997002+0.077375i  &  0.997002+0.077375i  &  0.997002+0.077375i  &  0.997002+0.077375i  &  0.997002+0.077375i  \\
          &       &       &       &  \\
    0.993511+0.113736i  &  0.993511+0.113736i  &  0.993511+0.113736i  &  0.993511+0.113736i  &  0.993511+0.113736i  \\
    0.993511+0.113736i  &  0.993511+0.113736i  &  0.993511+0.113736i  &  0.993511+0.113736i  &  0.993511+0.113736i  \\
    0.993511+0.113736i  &  0.993511+0.113736i  &  0.993511+0.113736i  &  0.993511+0.113736i  &  0.993511+0.113736i  \\
    0.993511+0.113736i  &  0.993511+0.113736i  &  0.993511+0.113736i  &  0.993511+0.113736i  &  0.993511+0.113736i  \\
    0.993511+0.113736i  &  0.993511+0.113736i  &  0.993511+0.113736i  &  0.993511+0.113736i  &  0.993511+0.113736i  \\
          &       &       &       &  \\
    0.988694+0.149945i  &  0.988694+0.149945i  &  0.988694+0.149945i  &  0.988694+0.149945i  &  0.988694+0.149945i  \\
    0.988694+0.149945i  &  0.988694+0.149945i  &  0.988694+0.149945i  &  0.988694+0.149945i  &  0.988694+0.149945i  \\
    0.988694+0.149945i  &  0.988694+0.149945i  &  0.988694+0.149945i  &  0.988694+0.149945i  &  0.988694+0.149945i  \\
    0.988694+0.149945i  &  0.988694+0.149945i  &  0.988694+0.149945i  &  0.988694+0.149945i  &  0.988694+0.149945i  \\
    0.988694+0.149945i  &  0.988694+0.149945i  &  0.988694+0.149945i  &  0.988694+0.149945i  &  0.988694+0.149945i  \\
    \bottomrule
    \end{tabular}%
   }
  \label{tab5}%
\end{table}%
\begin{table}[H]
  \centering
  \caption{Solution of (RED) at the grid points 5$\times$5$\times$5 in the cube.}
  {\fontsize{8.8}{9.6}\selectfont
    \begin{tabular}{rrrrr}
    \toprule
    0.999999+0.000010i  &  0.999999+0.000010i  &  0.999999+0.000010i  &  0.999999+0.000010i  &  0.999999+0.000010i  \\
    0.999999+0.000010i  &  0.999999+0.000010i  &  0.999999+0.000010i  &  0.999999+0.000010i  &  0.999999+0.000010i  \\
    0.999999+0.000010i  &  0.999999+0.000010i  &  0.999999+0.000010i  &  0.999999+0.000010i  &  0.999999+0.000010i  \\
    0.999999+0.000010i  &  0.999999+0.000010i  &  0.999999+0.000010i  &  0.999999+0.000010i  &  0.999999+0.000010i  \\
    0.999999+0.000010i  &  0.999999+0.000010i  &  0.999999+0.000010i  &  0.999999+0.000010i  &  0.999999+0.000010i  \\
          &       &       &       &  \\
    0.999332+0.036532i  &  0.999332+0.036532i  &  0.999332+0.036532i  &  0.999332+0.036532i  &  0.999332+0.036532i  \\
    0.999332+0.036532i  &  0.999332+0.036532i  &  0.999332+0.036532i  &  0.999332+0.036532i  &  0.999332+0.036532i  \\
    0.999332+0.036532i  &  0.999332+0.036532i  &  0.999332+0.036532i  &  0.999332+0.036532i  &  0.999332+0.036532i  \\
    0.999332+0.036532i  &  0.999332+0.036532i  &  0.999332+0.036532i  &  0.999332+0.036532i  &  0.999332+0.036532i  \\
    0.999332+0.036532i  &  0.999332+0.036532i  &  0.999332+0.036532i  &  0.999332+0.036532i  &  0.999332+0.036532i  \\
          &       &       &       &  \\
    0.997331+0.073005i  &  0.997331+0.073005i  &  0.997331+0.073005i  &  0.997331+0.073005i  &  0.997331+0.073005i  \\
    0.997331+0.073005i  &  0.997331+0.073005i  &  0.997331+0.073005i  &  0.997331+0.073005i  &  0.997331+0.073005i  \\
    0.997331+0.073005i  &  0.997331+0.073005i  &  0.997331+0.073005i  &  0.997331+0.073005i  &  0.997331+0.073005i  \\
    0.997331+0.073005i  &  0.997331+0.073005i  &  0.997331+0.073005i  &  0.997331+0.073005i  &  0.997331+0.073005i  \\
    0.997331+0.073005i  &  0.997331+0.073005i  &  0.997331+0.073005i  &  0.997331+0.073005i  &  0.997331+0.073005i  \\
          &       &       &       &  \\
    0.993999+0.109381i  &  0.993999+0.109381i  &  0.993999+0.109381i  &  0.993999+0.109381i  &  0.993999+0.109381i  \\
    0.993999+0.109381i  &  0.993999+0.109381i  &  0.993999+0.109381i  &  0.993999+0.109381i  &  0.993999+0.109381i  \\
    0.993999+0.109381i  &  0.993999+0.109381i  &  0.993999+0.109381i  &  0.993999+0.109381i  &  0.993999+0.109381i  \\
    0.993999+0.109381i  &  0.993999+0.109381i  &  0.993999+0.109381i  &  0.993999+0.109381i  &  0.993999+0.109381i  \\
    0.993999+0.109381i  &  0.993999+0.109381i  &  0.993999+0.109381i  &  0.993999+0.109381i  &  0.993999+0.109381i  \\
          &       &       &       &  \\
    0.989341+0.145611i  &  0.989341+0.145611i  &  0.989341+0.145611i  &  0.989341+0.145611i  &  0.989341+0.145611i  \\
    0.989341+0.145611i  &  0.989341+0.145611i  &  0.989341+0.145611i  &  0.989341+0.145611i  &  0.989341+0.145611i  \\
    0.989341+0.145611i  &  0.989341+0.145611i  &  0.989341+0.145611i  &  0.989341+0.145611i  &  0.989341+0.145611i  \\
    0.989341+0.145611i  &  0.989341+0.145611i  &  0.989341+0.145611i  &  0.989341+0.145611i  &  0.989341+0.145611i  \\
    0.989341+0.145611i  &  0.989341+0.145611i  &  0.989341+0.145611i  &  0.989341+0.145611i  &  0.989341+0.145611i  \\
    \bottomrule
    \end{tabular}%
   }
  \label{tab6}%
\end{table}%
\begin{table}[H]
  \centering
  \caption{Solution of (IE) at the grid points 5$\times$5$\times$5 in the cube.}
  {\fontsize{8.8}{9.6}\selectfont
    \begin{tabular}{rrrrr}
    \toprule
    1.000000+0.000010i  &  1.000000+0.000010i  &  1.000000+0.000010i  &  1.000000+0.000010i  &  1.000000+0.000010i  \\
    1.000000+0.000010i  &  1.000000+0.000010i  &  1.000000+0.000010i  &  1.000000+0.000010i  &  1.000000+0.000010i  \\
    1.000000+0.000010i  &  1.000000+0.000010i  &  1.000000+0.000010i  &  1.000000+0.000010i  &  1.000000+0.000010i  \\
    1.000000+0.000010i  &  1.000000+0.000010i  &  1.000000+0.000010i  &  1.000000+0.000010i  &  1.000000+0.000010i  \\
    1.000000+0.000010i  &  1.000000+0.000010i  &  1.000000+0.000010i  &  1.000000+0.000010i  &  1.000000+0.000010i  \\
          &       &       &       &  \\
    0.999332+0.036532i  &  0.999332+0.036532i  &  0.999332+0.036532i  &  0.999332+0.036532i  &  0.999332+0.036532i  \\
    0.999332+0.036532i  &  0.999332+0.036532i  &  0.999332+0.036532i  &  0.999332+0.036532i  &  0.999332+0.036532i  \\
    0.999332+0.036532i  &  0.999332+0.036532i  &  0.999332+0.036532i  &  0.999332+0.036532i  &  0.999332+0.036532i  \\
    0.999332+0.036532i  &  0.999332+0.036532i  &  0.999332+0.036532i  &  0.999332+0.036532i  &  0.999332+0.036532i  \\
    0.999332+0.036532i  &  0.999332+0.036532i  &  0.999332+0.036532i  &  0.999332+0.036532i  &  0.999332+0.036532i  \\
          &       &       &       &  \\
    0.997332+0.073005i  &  0.997332+0.073005i  &  0.997332+0.073005i  &  0.997332+0.073005i  &  0.997332+0.073005i  \\
    0.997332+0.073005i  &  0.997332+0.073005i  &  0.997332+0.073005i  &  0.997332+0.073005i  &  0.997332+0.073005i  \\
    0.997332+0.073005i  &  0.997332+0.073005i  &  0.997332+0.073005i  &  0.997332+0.073005i  &  0.997332+0.073005i  \\
    0.997332+0.073005i  &  0.997332+0.073005i  &  0.997332+0.073005i  &  0.997332+0.073005i  &  0.997332+0.073005i  \\
    0.997332+0.073005i  &  0.997332+0.073005i  &  0.997332+0.073005i  &  0.997332+0.073005i  &  0.997332+0.073005i  \\
          &       &       &       &  \\
    0.994000+0.109381i  &  0.994000+0.109381i  &  0.994000+0.109381i  &  0.994000+0.109381i  &  0.994000+0.109381i  \\
    0.994000+0.109381i  &  0.994000+0.109381i  &  0.994000+0.109381i  &  0.994000+0.109381i  &  0.994000+0.109381i  \\
    0.994000+0.109381i  &  0.994000+0.109381i  &  0.994000+0.109381i  &  0.994000+0.109381i  &  0.994000+0.109381i  \\
    0.994000+0.109381i  &  0.994000+0.109381i  &  0.994000+0.109381i  &  0.994000+0.109381i  &  0.994000+0.109381i  \\
    0.994000+0.109381i  &  0.994000+0.109381i  &  0.994000+0.109381i  &  0.994000+0.109381i  &  0.994000+0.109381i  \\
          &       &       &       &  \\
    0.989342+0.145611i  &  0.989342+0.145611i  &  0.989342+0.145611i  &  0.989342+0.145611i  &  0.989342+0.145611i  \\
    0.989342+0.145611i  &  0.989342+0.145611i  &  0.989342+0.145611i  &  0.989342+0.145611i  &  0.989342+0.145611i  \\
    0.989342+0.145611i  &  0.989342+0.145611i  &  0.989342+0.145611i  &  0.989342+0.145611i  &  0.989342+0.145611i  \\
    0.989342+0.145611i  &  0.989342+0.145611i  &  0.989342+0.145611i  &  0.989342+0.145611i  &  0.989342+0.145611i  \\
    0.989342+0.145611i  &  0.989342+0.145611i  &  0.989342+0.145611i  &  0.989342+0.145611i  &  0.989342+0.145611i  \\
    \bottomrule
    \end{tabular}%
   }
  \label{tab7}%
\end{table}%

\section*{Acknowledgements}
This work used the Extreme Science and Engineering Discovery Environment (XSEDE), which is supported by National Science Foundation grant number OCI-1053575.

\bibliographystyle{ieeetr}

\end{document}